\documentstyle[12pt]{article}
\begin{document}
\newtheorem{theorem}      {Th\'eor\`eme}[section]
\newtheorem{theorem*}     {theorem}
\newtheorem{proposition}  [theorem]{Proposition}
\newtheorem{definition}   [theorem]{Definition}
\newtheorem{e-lemme}        [theorem]{Lemma}
\newtheorem{cor}   [theorem]{Corollaire}
\newtheorem{resultat}     [theorem]{R\'esultat}
\newtheorem{eexercice}    [theorem]{Exercice}
\newtheorem{rrem}    [theorem]{Remarque}
\newtheorem{pprobleme}    [theorem]{Probl\`eme}
\newtheorem{eexemple}     [theorem]{Exemple}
\newcommand{\preuve}      {\paragraph{Preuve}}
\newenvironment{probleme} {\begin{pprobleme}\rm}{\end{pprobleme}}
\newenvironment{remarque} {\begin{rremarque}\rm}{\end{rremarque}}
\newenvironment{exercice} {\begin{eexercice}\rm}{\end{eexercice}}
\newenvironment{exemple}  {\begin{eexemple}\rm}{\end{eexemple}}
%
%
\newtheorem{e-theo}      [theorem]{Theorem}
\newtheorem{theo*}     [theorem]{Theorem}
\newtheorem{e-pro}  [theorem]{Proposition}
\newtheorem{e-def}   [theorem]{Definition}
\newtheorem{e-lem}        [theorem]{Lemma}
\newtheorem{e-cor}   [theorem]{Corollary}
\newtheorem{e-resultat}     [theorem]{Result}
\newtheorem{ex}    [theorem]{Exercise}
\newtheorem{e-rem}    [theorem]{Remark}
\newtheorem{prob}    [theorem]{Problem}
\newtheorem{example}     [theorem]{Example}
\newcommand{\proof}         {\paragraph{Proof~: }}
\newcommand{\hint}          {\paragraph{Hint}}
\newcommand{\heuristicproof}{\paragraph{heuristic proof}}
\newenvironment{e-probleme} {\begin{e-pprobleme}\rm}{\end{e-pprobleme}}
\newenvironment{e-remarque} {\begin{e-rremarque}\rm}{\end{e-rremarque}}
\newenvironment{e-exercice} {\begin{e-eexercice}\rm}{\end{e-eexercice}}
\newenvironment{e-exemple}  {\begin{e-eexemple}\rm}{\end{e-eexemple}}
\newcommand{\1}        {{\bf 1}}
\newcommand{\pp}       {{{\rm I\!\!\! P}}}
\newcommand{\qq}       {{{\rm I\!\!\! Q}}}
\newcommand{\B}        {{{\rm I\! B}}}
\newcommand{\cc}       {{{\rm I\!\!\! C}}}
\newcommand{\N}        {{{\rm I\! N}}}
\newcommand{\R}        {{{\rm I\! R}}}
\newcommand{\D}        {{{\rm I\! D}}}
\newcommand{\Z}       {{{\rm Z\!\! Z}}}
\newcommand{\C}        {{\bf C}}        
\newcommand{\CC}{{\cal C}}
\def \Re {{\rm Re\,}}
\def \Im {{ \rm Im\,}}
\def\Hom{{\rm Hom\,}}
\def\Lip{{\rm Lip}}
\def\st{{\rm st}}
\def\J{{\cal J}}
\def\A{{\cal A}}
\def\Mat{{\rm Mat}}
\def\m{{\rm m}}
\def\ind{{\rm ind\,}}
\def\<{\langle}
\def\>{\rangle}
\def\bar{\overline}
\def\Range{{\rm Range\,}}
\def\ker{{\rm ker\,}}
\def\Span{{\rm Span\,}}
\def\rank{{\rm rank\,}}
\def\id{{\rm id}}

%
%
\newcommand{\dontforget}[1]
{{\mbox{}\\\noindent\rule{1cm}{2mm}\hfill don't forget : #1
\hfill\rule{1cm}{2mm}}\typeout{---------- don't forget : #1 ------------}}
\newcommand{\note}[2]
{ \noindent{\sf #1 \hfill \today}

\noindent\mbox{}\hrulefill\mbox{}
\begin{quote}\begin{quote}\sf #2\end{quote}\end{quote}
\noindent\mbox{}\hrulefill\mbox{}
\vspace{1cm}
}
\title{Deformations and transversality \break
of pseudo-holomorphic discs}
\author{ Alexandre Sukhov{*} and Alexander Tumanov{**}}
\date{}
\maketitle

{\small
* Universit\'e des Sciences et Technologies de Lille, Laboratoire
Paul Painlev\'e,
U.F.R. de
Math\'e-matique, 59655 Villeneuve d'Ascq, Cedex, France,
 sukhov@math.univ-lille1.fr

** University of Illinois, Department of Mathematics
1409 West Green Street, Urbana, IL 61801, USA, tumanov@illinois.edu
}
\bigskip

Abstract.
We prove analogs of Thom's transversality theorem and
Whitney's theorem on immersions for pseudo-holomorphic discs.
We also prove that pseudo-holomorphic discs form a Banach
manifold.

\bigskip

MSC: 32H02, 53C15.

Key words: almost complex structure,  transversality, immersion,
jet, $J$-holomorphic  disc.
\bigskip

\section{Introduction}

Classical Thom's transversality theorem and Whitney's
approximation theorems (see, e. g., \cite{GG})
reflect the flexibility of smooth maps.
They are essentially local because global results
follow by means of cut-offs and the like.
In complex analytic category, in the absence of
cut-offs, the corresponding results
hold under global restrictions on the manifolds.
Kaliman and Zaidenberg \cite{KZ} prove that holomorphic
maps from a Stein manifold $X$ to a complex manifold $Y$
satisfy the jet transversality theorem after shrinking
the domain of the initial map.
Forstneri\v c \cite{Fo} proves that certain complex
analytic properties of the target manifold
(such as the existence of a dominating spray
in the sense of Gromov) imply global jet transversality
theorem (i.e., no shrinking of the domain is needed).
In almost complex setting,
the problem makes sense only when the source manifold
has complex dimension one, because holomorphic maps
of almost complex manifolds of higher dimension
generally do not exist.
Compact pseudo-holomorphic curves might not admit
any perturbations for a fixed almost complex structure.
In this paper we prove versions of the theorems of Thom
and Whitney for pseudo-holomorphic discs.
\begin{e-theo}
\label{Whitney}
Let $(M,J)$ be a $C^\infty$-smooth almost complex manifold
of complex dimension $n\ge2$ and
let $f_0:\D\to M$ be a $J$-holomorphic disc of class
$C^{m}(\bar\D)$, $m\ge0$; here $\D \subset \cc$ is the unit disc.
Then there exists a $J$-holomorphic immersion $f:\D\to M$
of class $C^\infty(\bar\D)$ arbitrarily close to $f_0$
in $C^{m}(\bar\D)$.
Furthermore, $f$ is homotopic to $f_0$ within
the set of $J$-holomorphic discs in $M$ close to $f_0$
in $C^{m}(\bar\D)$.
\end{e-theo}

This result is a version of Whitney's approximation theorem
for pseudo-holomorphic discs.
We need it in our work \cite{CoSuTu} on pseudo-holomorphic
discs in Stein domains.
We hope that Theorem \ref{Whitney} will find other
applications, in particular, in the theory of intersections
and moduli spaces of pseudo-holomorphic curves (see \cite{MS}).

For almost complex manifolds of complex dimension 2,
McDuff \cite{MD}
proves local approximation of a pseudo-holomorphic map by
an immersion near a singular point. The proof in \cite{MD} uses
a normal form of a holomorphic map at the singular point
(see also \cite{MicWh,Si}),
hence it uses the $C^\infty$ smoothness of $(M,J)$.
Although we state the result for $C^\infty$,
our proof goes through for finite smoothness.
McDuff \cite{MD}  also proves that a $J$-holomorphic map
$X\to(M,J)$ of a compact Riemann surface $X$ can be approximated
by a $J'$-holomorphic immersion for some $J'$ close to $J$.
In our Theorem \ref{Whitney}, the structure $J$ remains fixed.
Moreover, we do not separately analyze singular points
of the map.
Following the proof in the smooth category, we deduce
Theorem \ref{Whitney} from a pseudo-holomorphic
analog of Thom's transversality theorem below.

Let $D$ be an open subset in $\cc$ or
the closure of an open subset. By a $J$-complex
$k$-jet $D \to M$ at $\zeta \in D$ we mean the $k$-jet
of a smooth map $D \to M$ satisfying the equation of
$J$-holomorphicity (see equation (\ref{holomorphy}) below)
at $\zeta$ to order $k$. We denote by $J^k(D,M)$
the set of all $J$-complex $k$-jets $D \to M$.
Then $J^k(D,M)$ is a smooth submanifold of the space
of all $k$-jets of maps $D \to M$.

\begin{e-theo}
\label{Thom}
Let $(M,J)$ be a $C^\infty$-smooth  almost complex manifold.
Let $0\le k<m$ be integers.
Let $S$ be a $C^\infty$-smooth locally closed real manifold
in $J^k(\bar\D,M)$ with or without boundary.
Then the set of $J$-holomorphic
discs $f:\bar\D \to M$ of class $C^{m}(\overline\D)$
transverse to $S$
(that is, the $k$-jet extension $j^kf$ is transverse to $S$)
is dense in the space of all $J$-holomorphic discs
of class $C^{m}(\bar\D)$.
\end{e-theo}

We point out that $S$ is a real submanifold of $J^k(\bar\D,M)$
with no regard to the almost complex structures on $M$
or $J^k(\cc,M)$ (see \cite{LeS}).

Perturbations of small $J$-holomorphic discs are essentially
described by the implicit function theorem.
For big discs, Ivashkovich and Rosay \cite{IvRo}
prove that the center and direction of a $J$-holomorphic disc
admit arbitrary small perturbations.
Ivashkovich and Rosay \cite{IvRo} use the result in the study of
the Kobayashi-Royden metric on an almost complex manifold.

The proof of Theorem \ref{Thom} uses the Fredholm
property of the linearized Cauchy-Riemann operator
to parametrize $J$-holomorphic maps by holomorphic ones.
The main difficulty here is that in general the Fredholm
operator of the linearized problem has a non-trivial kernel.
We modify it  by adding a small holomorphic term
to obtain an invertible operator. The modification
is inspired by an important work of Bojarski \cite{Bo}.
The same idea yields another natural result.
\begin{e-theo}
\label{Manifold}
Let $(M,J)$ be a $C^\infty$-smooth  almost complex manifold.
The set of $J$-holomorphic discs of class $C^{k,\alpha}(\D)$,
$k\ge1$, $0<\alpha<1$ (resp. $W^{k,p}(\D)$, $k\ge1$, $2<p<\infty$)
forms a $C^\infty$-smooth Banach manifold modelled on the
space of holomorphic functions $\D\to\cc^n$ of the same class.
\end{e-theo}

For the usual complex structure,
Forstneri\v c \cite{Fo2} proved the result for holomorphic
mappings of a strongly pseudoconvex domain in a Stein manifold.
Finally, we point out (Proposition \ref{JonJdiscs})
that the manifold of
$J$-holomorphic discs carries a natural almost complex
structure only if $J$ is integrable.

Closing the introduction, we thank the referee for
useful critical remarks.

\section{Almost complex structures and
$J$-holomorphic discs}

We recall basic notions concerning almost complex manifolds
and pseudo-holomorphic discs.  Denote by $\D$  the
unit disc in $\cc$ and by $J_{\st}$  the standard complex structure
of $\cc^n$; the value of $n$ will be clear from the context.
Let $(M,J)$ be an almost complex manifold.
A smooth map $f:\D \to M$ is called {\it $J$-holomorphic}  if
\begin{eqnarray}
\label{CauchyRiemann}
 df \circ J_{\st} = J \circ df
\end{eqnarray}
We also call such a map  a $J$-{\it holomorphic} disc or
a {\it pseudo-holomorphic} disc.

In local coordinates $z\in\cc^n$, an almost complex structure
$J$ is represented by a $\R$-linear operator
$J(z):\cc^n\to\cc^n$, $z\in \cc^n$ such that $J(z)^2=-I$,
$I$ being the identity. Then the Cauchy-Riemann equations
(\ref{CauchyRiemann})  for a $J$-holomorphic disc $z:\D\to\cc^n$
can be written in the form
$$
z_\eta=J(z)z_\xi,\quad
\zeta=\xi+i\eta\in\D.
$$
We represent $J$ by a complex $n\times n$
matrix function $A=A(z)$ and obtain the equivalent equations
\begin{eqnarray}
\label{holomorphy}
z_{\bar\zeta}=A(z)\bar z_{\bar\zeta},\quad
\zeta\in\D.
\end{eqnarray}

We first recall the relation between $J$ and $A$ for
fixed $z$. Let $J:\cc^n\to\cc^n$ be a $\R$-linear map
so that $\det(J_\st+J)\ne0$, where $J_\st v=iv$.
Set
$$
Q=(J_\st+J)^{-1}(J_\st-J).
$$
One can show that $J^2=-I$ if and only if
$QJ_\st+J_\st Q=0$, that is, $Q$ is a complex anti-linear operator.
Then there is a unique matrix $A \in\Mat(n,\cc)$ such that
$$
Av=Q\bar v, \quad
v\in\cc^n.
$$
We introduce
$$
\J=\{J:\cc^n\to\cc^n: J\;{\rm is}\;\R{\rm-linear},\;J^2=-I,\;
\det(J_\st+J)\ne0  \},
$$
$$
\A=\{A\in\Mat(n,\cc): \det(I-A\bar A)\ne0 \}.
$$
It turns out that the map $J\mapsto A$ is a birational
homeomorphism $\J\to\A$ (see \cite{NiWo,SuTu}), and the
inverse map $A\mapsto J$ has the form
\begin{eqnarray}
\label{AtoJ}
Ju=i(I-A\bar A)^{-1}[(I+A\bar A)u-2A\bar u].
\end{eqnarray}
Let $J$ be an almost complex structure
in a domain $\Omega\subset\cc^n$.
Suppose $J(z)\in\J$, $z\in\Omega$. Then  $J$ defines
a unique complex matrix function $A$ in $\Omega$ such that
$A(z)\in\A$, $z\in\Omega$. We call $A$ the
{\it complex matrix} of $J$. The matrix $A$ has the
same regularity properties as $J$.
A function $f:\Omega\to\cc$ is $(J,J_\st)$-holomorphic
if and only if it satisfies the Cauchy-Riemann equations
\begin{eqnarray}
\label{CRscalar}
f_{\bar z}+f_z A=0,
\end{eqnarray}
where $f_{\bar z}$ and $f_z$ are considered row-vectors.
Such non-constant functions generally do not exist
unless $J$ is integrable; the integrability condition
in terms of $A=(a_{jk})$ has the form
\begin{eqnarray}
\label{integrabilityJ}
N_{jkl}=N_{jlk}, \quad
N_{jkl}:=(a_{jk})_{\bar z_l}+\sum_s (a_{jk})_{z_s}a_{sl}.
\end{eqnarray}

The main analytic tool in the study of pseudo-holomorphic discs
is the Cauchy-Green integral
$$
T u(\zeta) = \frac{1}{2\pi i}\int_{\D}
\frac{u(\omega)\,d\omega \wedge d\bar\omega}{\omega -\zeta}.
$$
As usual, we denote by $C^{k,\alpha}(\D)$ the space
of functions in $\D$ whose partial derivatives
to order $k$ satisfy a  H\"older condition with exponent
$0<\alpha<1$.
We write $C^\alpha(\D)=C^{0,\alpha}(\D)$.
We denote by $W^{k,p}(\D)$ the Sobolev space of functions
with derivatives to order $k$ in $L^p(\D)$.
We will use the following regularity properties of
the Cauchy-Green integral (see, e.g., \cite{Ve}).
\begin{proposition}
\label{operators}
\begin{itemize}
\item[(i)] Let $p > 2$ and $\alpha = (p-2)/p$. Then the linear
operator $T:L^p(\D)\to W^{1,p}(\D)$ is bounded.
The inclusion $W^{1,p}(\D)\subset C^\alpha(\D)$
is bounded, hence, the operator
$T:L^p(\D)\to C^\alpha(\D)$ is bounded, and
$T:L^p(\D)\to L^\infty(\D)$ is compact.
If $f\in L^p(\D)$, then $\partial_{\bar\zeta}Tf=f$,
$\zeta\in\D$, as a Sobolev derivative.
\item[(ii)] Let $1 < p < 2$ and $s= 2p(2-p)^{-1}$.
Then $T: L^p(\D) \to L^s(\D)$ is bounded.
\item[(iii)] Let $k \geq 0$
be integer and let $0<\alpha<1$. Then
$T:C^{k,\alpha}(\D) \to C^{k+1,\alpha}(\D)$
is bounded.
\end{itemize}
\end{proposition}

Using the Cauchy-Green operator we can replace the equations
(\ref{holomorphy}) by the equivalent integral equation
\begin{eqnarray}
\label{integralholomorphy}
z = T (A(z)\bar{z}_{\bar\zeta}) + \phi
\end{eqnarray}
where $\phi:\D \to \cc^n$ is an arbitrary holomorphic
vector function. In a small coordinate chart, we can
assume that the matrix $A$ is small with derivatives.
By the implicit function theorem, for every small
holomorphic vector function $\phi$, the equation
(\ref{integralholomorphy}) has a unique solution.
In particular, all $J$-holomorphic discs close
to a given small disc are parametrized by small
holomorphic vector functions.
We would like to establish a similar correspondence
for $J$-holomorphic discs which are not necessarily
small.

\section{Fredholm theory}

Let  $B_j$, $j = 1,2$  be $n\times n$ matrix functions
on $\D$ of class $L^p(\D)$, $p > 2$.
Solutions  of the equation
\begin{eqnarray}
\label{holvectors1}
u_{\bar\zeta} = B_1 u + B_2 \bar u
\end{eqnarray}
in the class $W^{1,p}(\D)$ are called
{\it  generalized holomorphic vectors}.
The equation (\ref{holvectors1}) arises as
the linearized equation (\ref{holomorphy}).
We introduce
$$T_0 u= Tu - Tu(0),$$
$$Pu = u - T_0(B_1 u + B_2 \bar u).$$
Let $r > 2p (p-2)^{-1}$. Then $s = (1/p + 1/r)^{-1} > 2$.
If $u \in L^r$, then $B_1 u, B_2 \bar u \in L^s$. Then
$T(B_1 u+ B_2 \bar u) \in W^{1,s}(\D)$ is continuous, and
$T_0(B_1 u+ B_2 \bar u)$ makes sense. Thus the operator
\begin{eqnarray}
\label{power}
P: L^r(\D) \to L^r(\D)
\end{eqnarray}
is bounded.
The equation (\ref{holvectors1}) is equivalent to the equation
\begin{eqnarray}
\label{holvectors5}
Pu = \phi
\end{eqnarray}
where $\phi$ is a usual $\cc^n$-valued holomorphic function
on the unit disc with $\phi(0) = u(0)$.
In the scalar case $n=1$ this equation admits a solution
for every $\phi$, and the kernel of $P$ is trivial.
This is a fundamental results of the theory of generalized analytic
functions \cite{Ve}. However, for $n > 1$ this is no longer true
(see \cite{Bo}). The equation (\ref{holvectors5}) in general
does not necessarily give a one-to-one correspondence between
generalized holomorphic vectors and usual holomorphic ones.

By Proposition \ref{operators}, the operator
$u \mapsto T_0(B_1 u + B_2 \bar u)$ is compact
and $P$ is Fredholm.
Hence the kernel of $P$ is finite-dimensional.
We modify the operator $P$ by adding a small holomorphic term
to obtain an operator with trivial kernel.
\begin{e-theo}
\label{centercontrol}
Let $B_j$, $j=1,2$ be $n \times n$ matrices of class
$L^p(\D)$, $p > 2$.
\begin{itemize}
\item[(i)]
Let $w_1,...,w_d$ form a basis of $\ker P$ over $\R$.
There exist holomorphic polynomial vectors
$p_1,...,p_d$ with $p_1(0)=...=p_d(0)=0$ such that  the operator
$\tilde P:L^r(\D)\to L^r(\D)$ of the form
\begin{eqnarray}
\label{modif}
\tilde P u = Pu + \sum_{j=1}^d (\Re (u,w_j))p_j
\end{eqnarray}
has trivial kernel. The polynomials $p_j$ can be chosen to be
arbitrarily small.
\item[(ii)]
If $B_j\in W^{k,p}(\D)$, $k\ge0$, $2<p<\infty$
(resp. $C^{k,\alpha}(\D)$, $k\ge0$, $0<\alpha<1$),
then $\tilde P$ is an invertible bounded operator
in the space $W^{k+1,p}(\D)$ (resp. $C^{k+1,\alpha}(\D)$),
in particular $\tilde P^{-1}$ is bounded.
The function $\phi=\tilde Pu$ is holomorphic if and only if
$u$ satisfies (\ref{holvectors1}). Furthermore,
$\tilde Pu(0)=u(0)$.
\end{itemize}
\end{e-theo}

We point out that a simpler version of Theorem \ref{centercontrol}
holds for the operator $T$ in place of $T_0$, but without the
conclusion $\tilde Pu(0)=u(0)$. We need the latter in the
proof of Theorem \ref{projets2} below.

In order to study the range of the operator $P$ we make use
of the adjoint $P^*$. For vector functions
$u = (u_1,...,u_n)$, $v = (v_1,...,v_n)$ we introduce
the usual inner product
$$
(u,v) = \sum_{j=1}^n\;
\frac{i}{2}\;
\int_\D u_j \bar v_j d\zeta \wedge d\bar\zeta.
$$
We consider the adjoints with respect to the real inner product
$\Re(\cdot,\cdot)$. The operator $P^*$ is defined on $L^q(\D)$
with $q = r (r-1)^{-1}$.
Given a function $\chi$ denote by $S_\chi$ the operator
$$
S_\chi u= \frac{1}{2\pi i}\int_\D \chi(\zeta) u(\zeta)
d\zeta \wedge d\bar\zeta
$$
We write $S_1$ for $\chi = 1$.
The following lemma is immediate.
\begin{e-lemme}
\begin{itemize}
\item[(i)] The adjoint of the operator of matrix multiplication
$u \mapsto B_j  u$ has the form  $u\mapsto B_j^*u$,
where $B_j^*$ denotes the hermitian transpose;
\item[(ii)] $T^* = -\bar T$,
here $\bar T u:=\bar{ T(\bar u)}$;
\item[(iii)] the conjugation operator
$\sigma: u \mapsto \bar u$ is self-adjoint: $\sigma^* = \sigma$;
\item[(iv)] $S_\chi^* = \bar\chi S_1$;
\item[(v)]  $T_0^* = -\bar T - \bar\zeta^{-1}S_1$;
and finally
\item[(vi)] $P^*=I+{\bar\zeta}^{-1}B_1^*\bar T\,\bar\zeta
+\zeta^{-1}\bar B_2^* T \zeta \sigma$.
\end{itemize}
\end{e-lemme}
\medskip

We need a boundary uniqueness theorem for solutions
of (\ref{holvectors1}), which in turn reduces to
the following version of the similarity principle
(\cite{GB}, Theorem 3.12; see also \cite{MS}, Theorem 2.3.5).
\begin{e-lemme}
\label{similarity}
Let $u\in W^{1,p}(\D)$ be a solution of (\ref{holvectors1}),
and let $\zeta_0\in b\D$. Then there exist a neighborhood
$U$ of $\zeta_0$ in $\cc$, a continuous nonsingular matrix function
$S$ in $\overline\D \cap U$,
and a vector function $\phi$ holomorphic in $\D \cap U$
and continuous in $\overline\D \cap U$ such that $u=S\phi$.
\end{e-lemme}

In \cite{GB}, $\zeta_0$ is an interior point,
but the proof goes through for a boundary point.
For completeness, we include a proof.
\proof
Put $B_1 = (a_{jk})$ and $B_2 = (b_{jk})$.
We rewrite the equation (\ref{holvectors1})
in the form
$$
u_{\bar\zeta} = H u.
$$
Here the matrix $H = (h_{jk})$ is defined
by $h_{jk} = a_{jk} + b_{jk}\bar{u_k}/u_k$,
and $u_k$ are the components of $u$.
Then $H \in L^p(\D)$.
For $\epsilon > 0$ put
$D(\zeta_0,\epsilon)=\{\zeta\in\bar\D :|\zeta-\zeta_0|
\le\epsilon \}$.
Define
$H_\epsilon(\zeta) = H(\zeta)$ for
$\zeta \in D(\zeta_0,\epsilon)$ and $H_\epsilon(\zeta) = 0$
otherwise.
Let $e_1,\ldots, e_n$ be the standard basis of  $\cc^n$.
Consider the equations
\begin{eqnarray}
\label{sim4}
w=T(H_\epsilon w) + e_j, \qquad j=1,\ldots, n.
\end{eqnarray}
By H\"older inequality (see \cite{Ve}, Theorem 1.19)
\begin{eqnarray}
\label{sim5}
||T(H_\epsilon w)||_{L^\infty(\D)}\le
C \epsilon^\alpha ||H||_{L^p(\D)} ||w||_{L^\infty(\D)},
\end{eqnarray}
where $\alpha=(p-2)/p$ and $C$ depends only on $p$.
Then for small $\epsilon>0$ the operator
$L:w \mapsto T(H_\epsilon w)$
is a contraction in $L^\infty(\D)$.
Then (\ref{sim4}) have unique solutions
$w_j\in L^\infty(\D)$.
By Proposition \ref{operators}(i), $w_j\in W^{1,p}(\D)$.
Define the matrix function $S = (w_1,...,w_n)$.
By (\ref{sim4}) and (\ref{sim5}),
$S\to I$ in $L^\infty(\D)$ as $\epsilon\to 0$;
here $I$ is the identity matrix.
Hence $S$ is nonsingular for small $\epsilon>0$.
Define $\phi = S^{-1} u$.
Then both $S$ and $\phi$ are in $W^{1,p}(\D)$,
in particular, they are continuous in $\bar\D$.
By differentiating (\ref{sim4}) we obtain
$S_{\bar\zeta}=HS$
in the interior of $D(\zeta_0,\epsilon)$.
We now have
$$
HS\phi=Hu=u_{\bar\zeta}=S_{\bar\zeta}\phi+S\phi_{\bar\zeta}
 =HS\phi+S\phi_{\bar\zeta}
$$
in the interior of $D(\zeta_0,\epsilon)$.
Since $S$ is nonsingular, then $\phi_{\bar\zeta} = 0$
in the interior of $D(\zeta_0,\epsilon)$.
The proof is complete.

\begin{e-cor}
\label{bunique}
Let $u\in W^{1,p}(\D)$ be a solution of (\ref{holvectors1}).
Suppose $u=0$ on an arc (or a set of positive length)
$E\subset b\D$. Then $u\equiv0$ identically.
\end{e-cor}
\proof
Let $\zeta_0\in b\D$ be a Lebesgue point of $E$.
Then $u$ has a decomposition $u=S\phi$ provided by
Lemma \ref{similarity}.
Since $S$ is nonsingular, then $\phi=0$ on $E$.
By the boundary uniqueness theorem
$\phi\equiv 0$, hence $u \equiv 0$ in a neighborhood of
$\zeta_0$.
By the connectedness argument, $u \equiv 0$ in all of $\D$,
as desired.
\medskip

Denote by $H_0$ the space of all holomorphic vector functions
$h\in L^r(\D)$ such that $h(0)=0$.
Then $H_0$ is a closed subspace of $L^r(\D)$.

\begin{e-lemme}
\label{range2}
$H_0+ \Range P = L^r(\D)$.
\end{e-lemme}
\proof
Let $v \in L^q(\D)$, $q = r (r-1)^{-1}$ be orthogonal
to both $H_0$ and $\Range P$. We prove that $v = 0$.

Since $(\Range P)^\bot=\ker P^*$, we have $P^*v = 0$, i.e.
\begin{eqnarray}
\label{boots}
\bar\zeta v = -B_1^* \bar{T(\zeta \bar v)}
-\bar\zeta \zeta^{-1}\bar B_2^* T(\zeta \bar v).
\end{eqnarray}
By ``bootstrapping'' we prove $\zeta \bar v \in L^p(\D)$.
Indeed, if $\zeta\bar v \in L^t(\D)$, $1 < t < 2$,
then $T(\zeta \bar v) \in L^s(\D)$
with $s = 2t (2-t)^{-1}$.
Then by (\ref{boots}) $\zeta \bar v \in L^k$ with
$$
k=\frac{1}{p^{-1}+s^{-1}}=\frac{t}{1-\beta t}
> \frac{t}{1 - \beta}
$$
where $\beta = \frac{1}{2} - \frac{1}{p}$,
$0 < \beta <\frac{1}{2}$.
Starting from $\zeta \bar v \in L^q(\D)$ and applying the
above argument  finitely many times,
we get $\zeta \bar v \in L^m(\D)$ for $m > 2$.
Finally, applying (\ref{boots}) one more time if necessary,
we get $\zeta \bar v \in L^p(\D)$.
Put $u = T(\zeta \bar v)$.
Then $u \in W^{1,p}(\D) \subset C^\alpha(\bar\D)$
with $\alpha = (p-2)p^{-1}$.

On the other hand, since $v$ is orthogonal to $H_0$,
then the function  $u = T(\zeta \bar v)$ vanishes
on $\cc \backslash \bar\D$ and, in particular, on $b\D$.
By (\ref{boots}), the function $u$ in $\D$ satisfies
\begin{eqnarray}
\label{sim1}
u_{\bar\zeta}=-\bar B_1^*u-\zeta \bar\zeta\,^{-1} B_2^* \bar u.
\end{eqnarray}
By Corollary \ref{bunique} applied to (\ref{sim1}),
we have $u \equiv 0$. Hence $v \equiv 0$ in $\D$.
The lemma is proved.
\medskip

{\bf Proof of Theorem \ref{centercontrol}:}
Part (i).
By Lemma \ref{range2}
there exist $p_1,...,p_d\in H_0$ such that
\begin{eqnarray}
\label{span}
\Span_{\R}(p_1,...,p_d) \oplus \Range P = L^r(\D)
\end{eqnarray}
By polynomial approximation, we can choose $p_j$ to be polynomial.
We now show that the operator $\tilde P$ defined by
(\ref{modif}) has trivial kernel.
Let $\tilde P u = 0$. Then by (\ref{span}) we have
$Pu = 0$ and $\Re (u,w_j) = 0$, $j= 1,...,d$.
Since the functions $w_1,...,w_d$ form a basis
of $\ker P$ over $\R$, we get $u = 0$. This
proves part (i). Part (ii) is immediate.
The theorem is proved.
\medskip

The following result is not new, and we do not need
it in the rest of the paper. We include it because
it is important by itself, and it is an immediate
consequence of Theorem \ref{centercontrol}.

\begin{e-cor}
Let $B_1,B_2\in L^p(\D)$, $p > 2$.
Then for every $\psi \in L^p(\D)$ the non-homogeneous equation
$$
u_{\bar\zeta} = B_1 u + B_2 \bar u + \psi
$$
has a solution in $W^{1,p}(\D)$.
\end{e-cor}
\proof
By Theorem 3.1, $u=\tilde P^{-1}T\psi\in W^{1,p}(\D)$
is a solution.

\section{Manifold of $J$-holomorphic discs}

Consider a non-homogeneous equation
\begin{eqnarray}
\label{CRg}
g_{\bar\zeta}=A(\zeta, g)\bar g_{\bar\zeta}+b(\zeta, g)
\end{eqnarray}
similar to the equation (\ref{holomorphy}).
In applications below, the matrix $A$ will arise from a certain
almost complex structure, hence we assume $\det(I-A\bar A)\ne0$.

\begin{proposition}
\label{Manifold-g}
Let $A$ and $b$ be $C^\infty$-smooth in an open set
in $\bar\D\times\cc^n$. Then the set of $C^{k,\alpha}(\D)$,
$k\ge1$, $0<\alpha<1$ (resp. $W^{k,p}(\D)$, $k\ge1$, $2<p<\infty$)
solutions of (\ref{CRg})
forms a $C^\infty$-smooth Banach manifold modelled on the
space of holomorphic functions $\D\to\cc^n$ of the same class.
\end{proposition}
\proof
Let $g_0$ be a solution of (\ref{CRg}).
We will construct a chart in a neighborhood of $g_0$.
Put $A_0(\zeta)=A(\zeta, g_0(\zeta))$. We make a substitution
by a real linear transformation
\begin{eqnarray}
\label{gtoh}
h=g-A_0\bar g.
\end{eqnarray}
The correspondence $g\leftrightarrow h$ is
one-to-one because
$$
g=(I-A_0\bar A_0)^{-1}(h+A_0\bar h).
$$
Then the equation (\ref{CRg}) turns into
\begin{eqnarray}
\label{CRh}
h_{\bar\zeta}=K_0\bar h_{\bar\zeta}+K_1 h+K_2\bar h+q.
\end{eqnarray}
The expressions of the new coefficients $K_0, K_1, K_2$, and $q$
in terms of $A$ and $b$ are obtained by straightforward computations.
In particular $K_0(\zeta, h_0(\zeta))=0$, which was
the goal of the substitution. Here $h_0=g_0-A_0\bar g_0$.

The equation (\ref{CRh}) is equivalent to
\begin{eqnarray}
\label{CRhT}
h=T_0(K_0\bar h_{\bar\zeta}+K_1 h+K_2\bar h+q)+\phi+\phi_0,
\end{eqnarray}
where $\phi$ is holomorphic, and $\phi_0$ is a fixed
holomorphic function such that (\ref{CRhT})
holds with $h=h_0$ and $\phi=0$.

We would like to apply the inverse function theorem
to the $C^\infty$ map
$$
h\mapsto F(h)=\phi=h-T_0(K_0\bar h_{\bar\zeta}
+K_1 h+K_2\bar h+q)-\phi_0
$$
in $C^{k,\alpha}(\D)$ (resp. $W^{k,p}(\D)$)
to obtain a one-to-one correspondence
$h\leftrightarrow \phi=F(h)$.

Note that $F(h_0)=0$.
The Fr\'echet derivative of the map $F$ at $h_0$ has the form
$$
F'(h_0)u=u-T_0(B_1u+B_2\bar u),
$$
where $B_1, B_2\in C^{k-1,\alpha}(\D)$ (resp. $W^{k-1,p}(\D)$).
Note that there is no term with $\bar u_{\bar\zeta}$ here
because $K_0(\zeta, h_0(\zeta))=0$.

In general the operator $F'(h_0)$  might not be an isomorphism
as required by the inverse function theorem.
We will modify the map $F$
using Theorem \ref{centercontrol}.

We apply Theorem \ref{centercontrol} to the operator
$$
P(u) =u-T_0(B_1 u+B_2\bar u),
$$
in $C^{k,\alpha}(\D)$ (resp. $W^{k,p}(\D)$)
to obtain an isomorphism
$$
\tilde P(u) =u-T_0(B_1 u+B_2\bar u)
+\sum_{j=1}^d\Re(w_j,u)p_j.
$$
We now modify the map $F$ (keeping the same notation $F$):
\begin{eqnarray}
\label{Fhtophi}
F(h)=\phi=h-T_0(K_0\bar h_{\bar\zeta}+K_1 h+K_2\bar h+q)-\phi_0
+\sum_{j=1}^d\Re(w_j,h-h_0) p_j.
\end{eqnarray}
Now $F'(h_0)=\tilde P$ is an isomorphism.
By the inverse function theorem $F^{-1}$ is well defined
and $C^{\infty}$-smooth in a neighborhood of zero in the space
of all vector functions of class $C^{k,\alpha}$
(resp. $W^{k,p}$). The map $F$ gives a one-to-one
correspondence between all solutions of (\ref{CRh})
close to $h_0$ and all holomorphic functions $\phi$
close to 0 in $C^{k,\alpha}$ (resp. $W^{k,p}$).
Hence it gives the desired chart.

The transition maps between the charts are smooth automatically
because the set of solutions of (\ref{CRg}) is a subset of
the manifold of all maps in each smoothness class, and the map $F$
is invertible as a map defined on all $C^{k,\alpha}$
(resp. $W^{k,p}$) functions close to $h_0$.

In conclusion we note that in this proof we could use
the operator $T$ instead of $T_0$ because a version
of Theorem \ref{centercontrol} holds for the operator $T$.
Proposition \ref{Manifold-g} is proved.
\medskip

{\bf Proof of Theorem \ref{Manifold}}.
Let $f_0:\D\to M$ be a $J$-holomorphic disc of class
$C^{k,\alpha}$ (resp. $W^{k,p}$).
We will construct a chart in a neighborhood of $f_0$.

Following \cite{IvRo},
by replacing each disc in $M$ by its graph in $\bar\D\times M$,
we reduce the proof to the case, in which the disc $f_0$
is an embedding and contained in a single chart in $\cc\times M$.
We further specify the choice of that chart.

For every point $p\in M$ there exists a chart
$\psi:U\subset M\to\cc^n$ such that $p\in U$, $\psi(p)=0$,
and for the push-forward
$\psi_*J = d\psi \circ J \circ d\psi^{-1}$ we have
$\psi_*J(0)=J_\st$.
We claim that we can choose such charts $\psi^\zeta$
for every point $p=f_0(\zeta)$ so that they
$C^{k,\alpha}$-smoothly depend on $\zeta\in\bar\D$.

Indeed, let $G:{\cal U}\subset T(M)\to M$ be a smooth map
defined in a neighborhood of the zero-section
of $T(M)$ such that for $G_p:=G|_{T_p(M)}:T_p(M)\to M$ we have
$G_{p\,*}(0)=\id$. For example, we can take for $G$
the exponential map for an arbitrary Riemannian
metric on $M$.
The pull-back bundle $f_0^*(T(M))\to\bar\D$ is trivial
as a complex vector bundle over $\bar\D$.
Then there exist $J$-complex linearly independent
$X_1(\zeta),\ldots, X_n(\zeta)\in T_p(M)$, $p=f_0(\zeta)$,
which are $C^{k,\alpha}$ in $\zeta\in\bar\D$
and $C^\infty$ in $\zeta\in\D$.
Define $\psi^\zeta(q)=z=(z_1,\ldots, z_n)\in\cc^n$
so that $G_p^{-1}(q)=\sum_{j=1}^n z_j X_j(\zeta)$,
$p=f_0(\zeta)$.
Here the product $z_j X_j$ is computed using $J(p)$.
Then $\psi^\zeta_*J(0)=J_\st$ because by construction
$\psi^\zeta_*$ is $(J,J_\st)$-linear at $p=f_0(\zeta)$.
The map  $\psi^\zeta$ is $C^{k,\alpha}$
in $\zeta\in\bar\D$ and $C^\infty$ in $\zeta\in\D$.

Using the map $(\zeta,q)\mapsto (\zeta,\psi^\zeta(q))$
we can introduce $C^{k,\alpha}$ coordinates in a neighborhood
of the graph of $f_0$ in $\bar\D\times M$. By shrinking
$\psi^\zeta(q)$ in $\zeta$ we obtain $C^\infty$ coordinates.

Let $0<r<1$.
Put $H(\zeta,q)=(\zeta,\psi^{r\zeta}(q))$,
here $\zeta\in\bar\D$ and $q\in M$ close to $f_0(\zeta)$.
For $r$ sufficiently close to 1, the $C^\infty$-smooth map
$$
H:U\subset\bar\D\times M\to\bar\D\times \cc^n
$$
defines coordinates $(\zeta,z)\in\bar\D\times\cc^n$ in
a neighborhood $U$ of the graph of $f_0$.

Consider the almost complex structure $J_\st\otimes J$ on
$\cc \times M$.  Then the push-forward
$\tilde J=H_*(J_\st\otimes J)$ is defined in a neighborhood
of $\bar\D\times \{0\}\subset\bar\D\times\cc^n$.
By the definition of $\tilde J$, we have
$\tilde J|_{\D\times \{0\}}=J_\st$, hence $\tilde J$
has complex matrix $\tilde A$ such that $\tilde A(\zeta,0)=0$.
Note that the projection $(\zeta,z)\mapsto \zeta$ is
$(\tilde J,J_\st)$-holomorphic, hence by (\ref{CRscalar})
the matrix $\tilde A$ has the form
$$
\tilde A = \left(
\begin{array}{cll}
0 &  0\\
b &  A
\end{array}
\right).
$$
The matrix $A(\zeta,\cdot)$ is the complex matrix
of the push-forward $\psi^{r\zeta}_*J$.
Since $\tilde A(\zeta,0)=0$, then
$A(\zeta,0)=0$ and $b(\zeta,0)=0$.
(The condition $\tilde A(\zeta,0)=0$ is not essential
in this proof, but it may be useful in other applications.)

A map $f:\D \to M$ close to $f_0$ is $J$-holomorphic
if and only if
$\zeta\mapsto(\zeta,f(\zeta))$
is $J_\st\otimes J$-holomorphic.
The latter is equivalent to
$\zeta\mapsto (\zeta,\psi^{r\zeta}(f(\zeta))$
being $\tilde J$-holomorphic.
Finally, it holds if and only if the map
$\zeta\mapsto g(\zeta):=\psi^{r\zeta}(f(\zeta))\in\cc^n$
satisfies the equation (\ref{CRg}), that is,
$$
g_{\bar\zeta}=A{\bar g}_{\bar\zeta}+b.
$$
Theorem \ref{Manifold} now follows by
Proposition \ref{Manifold-g}.

In conclusion we point out that if $k\ge2$, then
instead of using the substitution (\ref{gtoh}), we could pass
to the limit as $r\to 1$ and apply Theorem \ref{centercontrol}
in the limit.
The proof is complete.
\medskip

Let ${\cal M}$ be the manifold of all $J$-holomorphic discs in $M$
in a smoothness class admitted by Theorem \ref{Manifold}.
It seems reasonable to look for a suitable
almost complex structure on ${\cal M}$.
However, we show that a natural almost complex structure
exists on ${\cal M}$ only if $J$ is integrable.

Let $\Phi^\zeta: {\cal M}\ni f \mapsto f(\zeta)\in M$
be the evaluation map at $\zeta\in\D$.
\begin{e-pro}
\label{JonJdiscs}
Suppose there is an almost complex structure ${\cal J}$
on ${\cal M}$ such that the evaluation maps $\Phi^\zeta$
are $({\cal J},J)$-holomorphic. Then $J$ is integrable.
\end{e-pro}
\proof
Note that if ${\cal J}$ exists, then it is unique.
It suffices to prove the result for small discs
in a coordinate chart. Let $f\in{\cal M}$, that is
$f_{\bar\zeta}=A{\bar f}_{\bar\zeta}$.
Without loss of generality $f(0)=0$ and $A(0)=0$.
A tangent vector $u$ at $f$ is a solution of
the linearized equation, which implies
\begin{eqnarray}
\label{atzero}
u_{\bar\zeta}(0)=(A_z u+A_{\bar z}\bar u){\bar f}_{\bar\zeta}(0).
\end{eqnarray}
The hypothesis about ${\cal J}$ means that for every such $u$,
the map $Ju$ also satisfies that equation.
Using the description (\ref{AtoJ}) of $J$ in terms on $A$
and $A(0)=0$, we get
$Ju(0)=iu(0)$ and
$$
(Ju)_{\bar\zeta}(0)=i(u_{\bar\zeta}
-2A_z{\bar f}_{\bar\zeta}\bar u)(0).
$$
Applying (\ref{atzero}) to $Ju$ and comparing the result with
(\ref{atzero}) again we get
$$
A_{\bar z}{\bar f}_{\bar\zeta}\bar u(0)
{}=A_{\bar z}\bar u{\bar f}_{\bar\zeta}(0).
$$
Since $u(0)$ and ${\bar f}_{\bar\zeta}(0)$ are
arbitrary, we have
$$
(a_{jk})_{\bar z_l}(0)=(a_{jl})_{\bar z_k}(0).
$$
Since $A(0)=0$, then the integrability conditions
(\ref{integrabilityJ}) hold at 0, and the proposition
follows.

\section{Parametrization of jets and evaluation map}

We first prove the existence of solutions of (\ref{holvectors1})
with given jet at a fixed point.
\begin{e-theo}
\label{projets2}
Let $1\le k\le m$.
Let $B_1,B_2\in C^{m-1,\alpha}(\D)$, $0<\alpha<1$.
Then for every $a_0,...,a_k \in \cc^n$ there exists a solution $u$
of (\ref{holvectors1}) such that $u \in C^{m,\alpha}(\D)$ and
\begin{eqnarray}
\label{jets}
\partial_\zeta^j u(0) = a_j, \; 0 \leq j \leq k.
\end{eqnarray}
\end{e-theo}
\proof
Without loss of generality we assume $m=k$.
As a basis of induction, put $u=\tilde P^{-1}a_0$.
By Theorem \ref{centercontrol},
$u\in C^{k,\alpha}(\D)$ is a solution of (\ref{holvectors1})
with $u(0)=a_0$.
Assume by induction, that $u_1\in C^{k,\alpha}(\D)$
is a solution of (\ref{holvectors1}) with
\begin{eqnarray*}
\partial_\zeta^j u_1(0)=a_j,\; 0\leq j\le k-1.
\end{eqnarray*}
We look for a desired solution of (\ref{holvectors1}) in the form
\begin{eqnarray}
\label{jets4}
u = u_1 + \zeta^kv
\end{eqnarray}
Hence $v$ must satisfy the following conditions:
\begin{eqnarray}
\label{jets6}
v_{\bar\zeta} = B_1 v + B_2\zeta^{-k} \bar\zeta^k
 \bar v,
\end{eqnarray}
$$
v(0) = b,\quad
b = \frac{a_k - \partial_\zeta^k u_1(0)}{k!}.
$$
Let $\tilde P_0$ be the operator constructed
in Theorem \ref{centercontrol} for the equation (\ref{jets6}).
Define $v=\tilde P_0^{-1}b$.
We show that $u$ given by (\ref{jets4})
is in $C^{k,\alpha}(\D)$ and satisfies (\ref{jets}).

By Theorem \ref{centercontrol}, $v\in W^{1,p}(\D)$
for all $2<p<\infty$, hence $v\in C^\alpha(\D)$.
By bootstrapping, the equation $\tilde P_0v=b$ implies
$v\in C^{k,\alpha}(\D\setminus\D_{1/2})$.

Since $w=\zeta^kv$ satisfies (\ref{holvectors1}), then
\begin{eqnarray}
\label{jets66}
w = T(B_1 w + B_2\bar w)+\phi,
\end{eqnarray}
where $\phi\in C^\alpha(\D)$ is holomorphic in $\D$.

We claim that $\phi\in C^{k,\alpha}(b\D)$.
Indeed, we split the integral $T(B_1 w + B_2\bar w)$
over $\D$ into integrals over
$\D_{1/2}$ and $\D\setminus\D_{1/2}$.
The first term is holomorphic on $b\D$.
The second term is in $C^{k,\alpha}(\D\setminus\D_{1/2})$
by Proposition \ref{operators}(iii)
because $w=\zeta^kv\in C^{k,\alpha}(\D\setminus\D_{1/2})$,
and $B_1, B_2 \in C^{k-1,\alpha}(\D)$.
Hence the claim follows.

Since $\phi$ is holomorphic in $\D$ and $\phi\in C^{k,\alpha}(b\D)$,
then by the regularity of the Poisson integral,
$\phi\in C^{k,\alpha}(\D)$.
Finally by bootstrapping, using again
Proposition \ref{operators}(iii),
the equation (\ref{jets66})
yields $w\in C^{k,\alpha}(\D)$.
Hence $u\in C^{k,\alpha}(\D)$.

Since $v \in C^\alpha(\D)$, then Taylor's expansion of $w$
at 0 has the form
$$
w(\zeta)=b\zeta^k+O(|\zeta|^{k+\alpha}).
$$
Hence, $\partial_\zeta^k w(0) = k!\,b$
and $\partial_\zeta^j w(0)=0$ for $0\le j<k$.
Hence, $u = u_1  + \zeta^kv$ satisfies (\ref{jets}).
The proof is complete.
\bigskip

Let $j^k_\zeta: (C^k(\bar\D))^n \to (\cc^n)^{k+1}$ be
holomorphic $k$-jet evaluation map at $\zeta \in \bar\D$
defined by
$$
j_\zeta^k(u) = \{ \partial_\zeta^ju(\zeta): 0 \leq j \leq k \}.
$$
Let ${\cal V }^{m,\alpha}$ denote the space of
all solutions of (\ref{holvectors1}) of class $C^{m,\alpha}(\D)$.

\begin{e-pro}
\label{evaluation}
Let $1\le k \le m$.
Let $B_1, B_2 \in C^{m-1,\alpha}(\D)$, $0 < \alpha < 1$.
Then there exists a subspace $V \subset {\cal V}^{m,\alpha}$,
$\dim V < \infty$,
such that for every $\zeta \in \bar\D$ the restriction
$j^k_\zeta|_V: V \to (\cc^n)^{k+1}$ is surjective.
\end{e-pro}

\proof
By Theorem \ref{projets2} there exists a subspace
$V_0 \subset{\cal V}^{m,\alpha}$ such that $j_0^k\vert_{V_0}$
is bijective.
The evaluation point $0$ in Theorem \ref{projets2} can be replaced
by any other point $\zeta \in \bar\D$ by extending $B_j$
to all of $\cc$ and applying Theorem \ref{projets2} to a bigger
disc with center at $\zeta$.
Let $V(\zeta) \subset {\cal V}^{m,\alpha}$
denote such a subspace that $j^k_\zeta\vert_{V(\zeta)}$ is bijective.
By continuity, every $\zeta_0 \in \bar\D$ has a neighborhood
$U(\zeta_0) \subset \D$ so that for every $\zeta \in U(\zeta_0)$
the map $j^k_\zeta|_{V(\zeta_0)}$ is bijective. By compactness,
we find a finite covering $\cup_{l=1}^N U(\zeta_l)$ of $\bar\D$
and then $V = \sum_{l=1}^N V(\zeta_l)$ has the needed property.
The proof is complete.

\section{Transversality of $J$-holomorphic discs}

In this section we prove Theorems \ref{Whitney} and \ref{Thom}.

Let $(M,J)$ be an almost complex manifold.
Recall that by a $J$-complex $k$-jet $\bar\D \to M$
at $\zeta\in\bar\D$ we mean the $k$-jet of a smooth map
$\bar\D \to M$ $J$-holomorphic to order $k$ at $\zeta$.
We denote by $J^k(\bar\D,M)$ the manifold of $J$-complex $k$-jets
$\bar\D \to M$.
Let $U \subset M \to \cc^n$ be a coordinate chart in which $J$
has complex matrix $A$.
Let $f: \bar\D \to U$ be a smooth map $J$-holomorphic to
order $k$ at $\zeta$.
Then all the derivatives
$\partial^p_\zeta\partial^q_{\bar\zeta}f(\zeta)$,
$p + q \leq k$, are uniquely determined by
$\partial_\zeta^p f(\zeta)$, $p \leq k$ from the equation
(\ref{holomorphy}). Hence for every chart $U \subset M$,
in which $J$ is defined by its complex matrix $A$, and
open set $D \subset \bar\D$, there is a chart
$\tilde U \subset J^k(\bar\D,M) \to D \times U \times\cc^{nk}$
so that the $k$-jet extension of a $J$-holomorphic map
$f:D \to U$ has the representaion
$$
j^k f(\zeta)=(\zeta, f(\zeta), \partial_\zeta^1f(\zeta),
 \ldots, \partial_\zeta^kf(\zeta)).
$$

{\bf Proof of Theorem \ref{Thom}.}
Let $f_0:\D\to M$ be a $J$-holomorphic disc of class
$C^{m}(\bar\D)$, $m>k$.
We would like to find a $J$-holomorphic disc $f:\D\to M$
close to $f_0$ in $C^{m}(\bar\D)$ so that $j^kf$
is transverse to $S$.

Without loss of generality we can assume that $m$ is large.
Indeed, for $0<r<1$ the map $\zeta\mapsto f_0(r\zeta)$
is in $C^\infty(\bar\D)$ and approaches $f_0$ in
$C^{m}(\bar\D)$ as $r\to 1$.

We will use elements of the constructions from the proofs
of Theorem \ref{Manifold} and Proposition \ref{Manifold-g}.
In particular, we will use the one-to-one correspondence
$f\leftrightarrow h$ between $J$-holomorphic discs $f$
close to $f_0$ and solutions $h$ of the equation (\ref{CRh})
close to $h_0$.
We will also use the map $F:h\mapsto\phi$ defined by
(\ref{Fhtophi}).
For this reason, we assume that $f_0\in C^{m,\alpha}(\D)$,
where $m$ is large and $0<\alpha<1$. The particular values
of $m$ and $\alpha$ are unimportant.

Similarly to $J^k(\bar\D,M)$,
we introduce $\tilde J^k(\bar\D,\cc^n)$ as the manifold of
$k$-jets of smooth maps satisfying (\ref{CRh})
to order $k$ at a point.
Then $\tilde J^k(\bar\D,\cc^n)$ consists of a single chart
$\tilde J^k(\bar\D,\cc^n)\ni j^k h(\zeta)\mapsto
(\zeta, h(\zeta), \partial_\zeta^1 h(\zeta),
\ldots, \partial_\zeta^k h(\zeta))\in\cc^{n(k+1)+1}$.

The correspondence $f\leftrightarrow h$ gives rise to
a $C^\infty$-diffeomorphism
$$
\Psi:J^k(\bar\D,M)\supset U\to
\tilde U\subset\tilde J^k(\bar\D,\cc^n)
$$
defined in a neighborhood $U\supset j^k f_0(\bar\D)$.
Then $j^kf$ is transverse to $S$ in $J^k(\bar\D,M)$
if and only if $j^k h$ is transverse to $\tilde S:=\Psi(S)$
in $\tilde J^k(\bar\D,\cc^n)$.

Let $V$ be the space of solutions of the equation
(\ref{holvectors1})
provided by Proposition \ref{evaluation}.
Let $u_1,...,u_N$ form a basis of $V$.
Let $\phi_j =F'(h_0)(u_j)$.
For $s = (s_1,...,s_N) \in \R^N$ put
$\phi_s =\sum_l s_l \phi_l$. Then the map
\begin{eqnarray*}
& &\Phi:\bar\D \times \R^N \to \tilde J^k(\bar\D,\cc^n),\\
& &\Phi(\zeta,s) = j^kF^{-1}(\phi_s)(\zeta)
\end{eqnarray*}
is defined for small $s \in \R^N$.
The map $\Phi$ is $C^{m-k,\alpha}$-smooth.
By Proposition \ref{evaluation},
the map $\Phi$ is a submersion for $s = 0$, $\zeta \in \bar\D$.
Then $X=\Phi^{-1}(\tilde S)$ is a $C^{m-k}$-smooth manifold
in a neighborhood of
$\bar\D \times \{0 \} \subset \bar\D \times \R^N$.
Let $\tau: \bar\D \times \R^N\to \R^N$ be the projection.
Then $\Phi(\cdot,s)$ is transverse
to $\tilde S \subset\tilde J^k(\bar\D,\cc^n)$
if and only if $s$ is a regular value of $\tau|_X$.
(See the proof of Thom's transversality theorem
in \cite{GG}, Lemma 4.6, Chapter 2.)

Since $m$ is large,
then by Sard's theorem the set of critical values
of $\tau|_X$ has measure zero.
Hence there exists $s\in\R^N$ arbitrarily close to $0$
such that $h=F^{-1}(\phi_s)$ has $k$-jet extension
$\Phi(\cdot,s)$ transverse
to $\tilde S$ in $\tilde J^k(\bar\D,\cc^n)$.
Then the corresponding $j^kf$ is transverse to $S$ in
$J^k(\bar\D,M)$, as desired.
Theorem \ref{Thom} is now proved.
\medskip

{\bf Proof of Theorem \ref{Whitney}.}
Theorem \ref{Whitney} follows from Theorem \ref{Thom}
in the same way that the classical Whitney theorem
follows from that of Thom (see \cite{GG}).

The space $J^1(\bar\D,M)$ has a canonical structure
of a vector bundle over $\bar\D\times M$.
Let $S$ be the zero section of this bundle.
Then a $J$-holomorphic map $f:\D \to M$ is an immersion
if and only if $j^1f(\bar\D)\cap S=\emptyset$.
We have $\dim_\R J^1(\bar\D,M) = 4n + 2$ and
$\dim_\R S= 2n + 2$.
If $j^1f$ is transverse to $S$, then indeed
$j^1f(\bar\D)\cap S=\emptyset$, because
otherwise $\dim_\R \D+\dim_\R S\ge\dim_\R J^1(\bar\D,M)$
implies $n\le1$.
The desired result now follows by Theorem \ref{Thom}.
The homotopy statement follows by Theorem \ref{Manifold}.
The proof is complete.

\end{document}